\documentclass[12pt]{amsart}
\usepackage{amscd}
\usepackage{graphicx}
\def\frk{\frak}               % font for "Fraktur"

\def\Phi{{\frk n}}
\def\Phi{{\frk N}}
\def\opn#1#2{\def#1{\operatorname{#2}}} % to make operators
\opn\chara{char} \opn\length{\ell} \opn\pd{pd} \opn\rk{rk}
\opn\projdim{proj\,dim} \opn\injdim{inj\,dim} \opn\rank{rank}
\opn\depth{depth} \opn\grade{grade} \opn\height{height}
\opn\embdim{emb\,dim} \opn\codim{codim}

\opn\Tr{Tr} \opn\bigrank{big\,rank}
\opn\superheight{superheight}\opn\lcm{lcm}
\opn\trdeg{tr\,deg}%\emph{
\opn\reg{reg} \opn\lreg{lreg} \opn\ini{in} \opn\lpd{lpd}
\opn\size{size}
\opn\div{div} \opn\Div{Div} \opn\cl{cl} \opn\Cl{Cl}
\opn\Spec{Spec} \opn\Supp{Supp} \opn\supp{supp} \opn\Sing{Sing}
\opn\Ass{Ass} \opn\Min{Min}
\opn\Ann{Ann} \opn\Rad{Rad} \opn\Soc{Soc}
\opn\Im{Im} \opn\Ker{Ker} \opn\Coker{Coker} \opn\Am{Am}
\opn\Hom{Hom} \opn\Tor{Tor} \opn\Ext{Ext} \opn\End{End}
\opn\Aut{Aut} \opn\id{id}

\opn\nat{nat}
\opn\pff{pf}%   \pf exists already
\opn\Pf{Pf} \opn\GL{GL} \opn\SL{SL} \opn\mod{mod} \opn\ord{ord}
\opn\Gin{Gin} \opn\Hilb{Hilb}
\opn\aff{aff} \opn\con{conv} \opn\relint{relint} \opn\st{st}
\opn\lk{lk} \opn\cn{cn} \opn\core{core} \opn\vol{vol}
\opn\link{link} \opn\star{star}
\opn\gr{gr}

\def\pot#1#2{#1[\kern-0.28ex[#2]\kern-0.28ex]}

\opn\dirlim{\underrightarrow{\lim}}
\opn\inivlim{\underleftarrow{\lim}}
\def\Implies{\ifmmode\Longrightarrow \else
        \unskip${}\Longrightarrow{}$\ignorespaces\fi}
\def\implies{\ifmmode\Rightarrow \else
        \unskip${}\Rightarrow{}$\ignorespaces\fi}
\def\iff{\ifmmode\Longleftrightarrow \else
        \unskip${}\Longleftrightarrow{}$\ignorespaces\fi}

\let\:=\colon
\newtheorem{Theorem}{Theorem}[section]
\newtheorem{Lemma}[Theorem]{Lemma}

\newtheorem{Remark}[Theorem]{Remark}

\newtheorem{Example}[Theorem]{Example}

\newtheorem{Definition}[Theorem]{Definition}

\let\epsilon\varepsilon
\let\phi=\varphi
\let\kappa=\varkappa
\def\qed{\ifhmode\textqed\fi
      \ifmmode\ifinner\quad\qedsymbol\else\dispqed\fi\fi}
\def\textqed{\unskip\nobreak\penalty50
       \hskip2em\hbox{}\nobreak\hfil\qedsymbol
       \parfillskip=0pt \finalhyphendemerits=0}
\def\dispqed{\rlap{\qquad\qedsymbol}}

\opn\dis{dis}
\def\pnt{{\raise0.5mm\hbox{\large\bf.}}}

\opn\Lex{Lex}

\begin{document}

\title{$f$-ideals of degree $2$}

\author{G. Q. Abbasi*, S. Ahmad*, I. Anwar*, W. A. Baig*}

%\thanks{The author is highly grateful to the Abdus Salam School of Mathematical Sciences,
 % GC University, Lahore, Pakistan in supporting and facilitating this
  %research. The author  would like to thank   Prof. J\"{u}rgen Herzog for introducing the
  %idea and encouragement.}

\address{* COMSATS Institute of Information Technology,
    Lahore, Pakistan.}
\email { qanberabbasi@ciitlahore.edu.pk, sarfraz11@gmail.com,
     iimrananwar@gmail.com, waqqasbaig82@gmail.com.}
 \maketitle

\begin{abstract}
In this paper, we introduce the concept of $f$-ideals and discuss
its algebraic properties. In particular, we give the
characterization of all the $f$-ideals of degree $2$.
% Also for a simplicial complex
%$\Delta$, we have discussed totally link vector $\l(\Delta)$ and its
%properties.
 \vskip 0.4 true cm
 \noindent
\end{abstract}

\section{Introduction}

%A finite simplicial complex $\Delta$ consists of a finite set of
%vertices $V = \{x_1, . . . ,x_n\}$ is a collection of subsets of $V$
%satisfying the following axioms;\\
%$(i)$ $\{x_i\}\in \Delta$ for all $i,\, \, \, 1\leq i\leq n$.\\
%$(ii)$ If $F\in \Delta$ and $H\subset F$, then $H\in \Delta$.\\

% For a simplicial complex
%$\Delta$, combinatorial properties under consideration are
%$f$-vector.
% and the total link vector
 %$\l(\Delta)$.
  %Whereas the
%algebraic properties discussed in this paper are mainly dealt with
%the characterization of $f$-ideals of degree $2$.\\
The aim of this paper is to explore the algebraic and combinatorial
properties of simplicial complexes. Let $S=k[x_1, . . . ,x_n]$ be a
polynomial ring over an infinite field $k$.  There is a natural
bijection between a square-free monomial ideal and a simplicial
complex:
$$\Delta \leftrightarrow I_{\mathcal{N}}$$
Where $I_{\mathcal{N}}$ is known as the Stanley Reisner ideal or
non-face ideal of $\Delta$. This one to one correspondence has been
discussed widely in the literature for instance in
\cite{BH},\cite{Mi}, \cite{St} and
\cite{Vi}.\\
%Stanley Reisner ideal $I_{\mathcal{N}}$ gives rise to the
%Stanley Reisner ring an algebraic object related to the simplicial
%complex $\Delta$ introduced by R.P. Stanley in \cite{St1}.\\

In \cite{F1}, Faridi introduced another correspondence:
$$\Delta \leftrightarrow I_{\mathcal{F}}$$
Where $I_{\mathcal{F}}$ is the facet ideal of a given simplicial
complex $\Delta$. In \cite{F1} and \cite{F2} Faridi has discussed
and investigated some algebraic compatibilities of these two ideals
for a given simplicial complex $\Delta$. Also when the simplicial
complex is a tree (defined in \cite{F1}), its facet ideal posseses
interesting algebraic
and combinatorial properties discussed in \cite{F1} and \cite{F2}. \\

Given a square free monomial ideal $I$, one can consider it as the
facet ideal of a simplicial complex $\delta_{\mathcal{F}}(I)$, and
the Stanley-Reisner ideal of another $\delta_{\mathcal{N}}(I)$. So
for a square free monomial ideal $I$, one can explore some
invariants of $\delta_{\mathcal{F}}(I)$ and
$\delta_{\mathcal{N}}(I)$.\\

 In this paper, We introduce the
$f$-ideals and in Theorem \ref{Th} we give the characterization of
all the $f$-ideals of degree $2$.  A monomial ideal is called
$f$-ideal if both the simplicial complexes $\delta_{\mathcal{F}}(I)$
and $\delta_{\mathcal{N}}(I)$ have the same $f$-vector, where
$f$-vector of a $d$ dimensional simplicial complex $\Delta$ is the
$(d+1)$-tuple:
$$f(\Delta)=(f_0 , f_1 ,\ldots , f_d),$$
where $f_i$ is the number of faces of dimension $i$ of $\Delta$.

\section{Basic combinatorics and algebra of simplicial complexes}
This section is a review on the combinatorics and algebra associated
to simplicial complexes discussed in [\cite{F1} - \cite{Mi} and
\cite{Vi}].
\begin{Definition}
\em{A  simplicial complex $\Delta$ over a set of vertices
$V=\{x_1,x_2,\ldots,x_n\}$ is a collection of subsets of $V$, with
the property that $\{x_i\}\in \Delta$ for all $i$, and if $F\in
\Delta$ then all subsets of $F$ are also in $\Delta$(including the
empty set). An element of $\Delta$ is called a {\em face} of
$\Delta$, and the {\em dimension of a face} $F$ of $\Delta$ is
defined as $|F|-1$, where $|F|$ is the number of vertices of $F$.
The faces of dimension $0$ and $1$ are called {\em vertices and
edges}, respectively, and $dim\emptyset=-1$. The maximal faces of
$\Delta$ under inclusion are called {\em facets}. }
\end{Definition}
We denote simplicial complex $\Delta$ by a generating set of its
facets $F_1,\ldots  , F_q $ as
$$\Delta = <F_1,\ldots , F_q>$$
Also, we denote the {\em facet set} by $\mathcal{F}=\{F_1,\ldots ,
F_q\}$. A simplicial complex with only one facet is called a {\em
simplex}.

The following definitions lay the foundation of the dictionary
between the combinatorial and algebraic properties of the simplicial
complexes over the finite set of vertices $[n]$.

\begin{Definition}
{\em Let $\Delta$ be a simplicial complex over $n$ vertices
$\{v_1,\ldots , v_n\}$. Let $k$ be a field, $x_1, . . . , x_n$ be
indeterminates, and $S$ be the polynomial ring
$k[x_1,\ldots, x_n]$. \\
(a) We define $I_{\mathcal{F}}$ to be the ideal of $S$ generated by
square-free monomials $x_{i1} \ldots x_{is}$, where $\{v_{i1} ,
\ldots, v_{is}\}$ is a facet of $\Delta$. We call
$I_{\mathcal{F}}$ the facet ideal of $\Delta$.\\
(b) We define $I_{\mathcal{N}}$ to be the ideal of $S$ generated by
square-free monomials $x_{i1} \ldots x_{is}$ , where $\{v_{i1} ,
\ldots , v_{is}\}$ is not a face of $\Delta$. We call
$I_{\mathcal{N}}$ the non-face ideal or the Stanley- Reisner ideal
of $\Delta$. }
\end{Definition}

\begin{Definition}\label{fc}
{\em Let $I = (M_1,\ldots,M_q)$ be an ideal in a polynomial ring
$k[x_1, . . . , x_n]$, where $k$ is a field and $M_1,\ldots ,M_q$
are square-free monomials in $x_1,\ldots , x_n$ that form a minimal
set
of generators for $I$.\\
(a) We define $\delta_{\mathcal{F}}(I)$ to be the simplicial complex
over a set of vertices $v_1,\ldots , v_n$ with facets $F_1,\ldots ,
F_q$, where for each $i$, $F_i = \{v_j \, \, |\, \,  x_j |M_i, 1
\leq j\leq n\}$. We call
$\delta_{\mathcal{F}}(I)$ the {\em facet complex of $I$.}\\
(b) We define $\delta_{\mathcal{N}}(I)$ to be the simplicial complex
over a set of vertices $v_1,\ldots , v_n$, where $\{v_{i1},\ldots ,
v_{is}\}$ is a face of $\delta_{\mathcal{N}}(I)$ if and only if
$x_{i1} . . . x_{is} \not \in I$. We call $\delta_{\mathcal{N}}(I)$
the {\em non-face complex} or {\em the Stanley-Reisner complex} of
$I$. }
\end{Definition}
\begin{Remark}{\em
For given a square free monomial ideal $I$, one can construct
$\delta_{\mathcal{F}}(I)$ by using the above definition (a). Where
Faridi in \cite{F1}, has given the construction of
$\delta_\mathcal{N}(I)$ by using the minimal vertex cover of
$\delta_\mathcal{F}(I)$.}
\end{Remark}
\begin{Example}\label{ex}

Let $I=(xy,yz)\subset k[x,y,z]$, then following are non-face complex
and facet complex.

\begin{figure}[h]
\begin{center}
  \includegraphics[width=8cm]{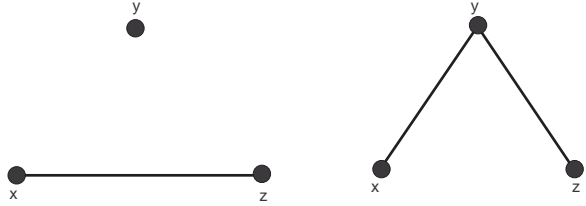}\vspace{0.0cm}
  \caption{Non-face and facet complex}\label{}
\end{center}
\end{figure}

\end{Example}
%One may think about the following proposition to construct
%$\delta_{\mathcal{N}}(I)$ as well:
%\begin{Proposition}\label{ps}\em{ [\cite{Vi},Prop 5.3.10]}

%If $\Delta$ is a simplicial complex with vertices $x_1,\ldots,x_n$,
%then the primary decomposition of the Stanley-Reisner ideal of
%$I_{\mathcal{N}}$ is:
%$$I_{\mathcal{N}}=\bigcap_{F}P_F$$
%where the intersection is taken over all facets $F$ of $\Delta$, and
%$P_\mathcal{F}$ denotes the face ideal generated by all $x_i$ such
%that $x_i\not\in F$.

%\end{Proposition}

\begin{Definition} {\em Let $S=k[x_1,\ldots,x_n]$ be a polynomial ring,
the {\em Support} of a monomial $x^a =x_1^{a_1}\ldots x_n^{a_n}$ in
$S$ is given by $Supp(x^a)=\{ x_i|a_i>0\}$.\\
Similarly, let $I=(g_1,\ldots,g_m)\subset S$ be a square-free
monomial ideal then $$Supp(I)=\bigcup_{i=1}^{m}Supp(g_i)$$ }
\end{Definition}
\begin{Remark}{\em
It is worth noting that for any square-free monomial ideal $I\subset
S$ the $\delta_{\mathcal{F}}(I)$ will be a simplicial complex on the
vertex set $[s]$, where $s=|Supp(I)|$. But $\delta_{\mathcal{N}}(I)$
will be a simplicial complex on $[n]$. So both
$\delta_{\mathcal{F}}(I)$ and $\delta_{\mathcal{N}}(I)$ will have
the same vertex set if and only if $Supp(I)=\{x_1,\ldots,x_n\}$.\\
For example,  for the ideal $I=(x_2x_3, x_2x_4, x_3x_4)$ in
$S=k[x_1,x_2,x_3,x_4]$, $$\delta_{\mathcal{F}}(I)=<\{v_2,v_3\},
\{v_2,v_4\}, \{v_3,v_4\}>$$ is a simplicial complex on the vertex
set $\{v_2,v_3,v_4\}$. Whereas,
$$\delta_{\mathcal{N}}(I)=<\{v_1,v_4\},\{v_1,v_3\}\{v_1,v_2\}>$$ is
the simplicial complex on the vertex set $\{v_1, v_2,v_3,v_4\}$.

 \begin{figure}[h]
\begin{center}
  \includegraphics[width=8cm]{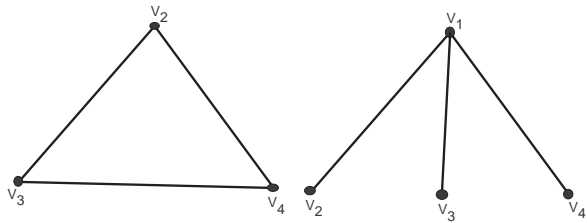}\vspace{0.0cm}
  \caption{facet and non-face complex}\label{}
\end{center}
\end{figure}}
\end{Remark}
\begin{Definition}{\em
Let $I=(g_1,\ldots,g_m)\subset S$ be a square-free monomial ideal,
the $deg(I)$ is defined as:
$$deg(I)=Sup\{deg(g_i)|i\in\{1,\ldots,m\}\}$$}
\end{Definition}
\section{$f$-ideals and classification of $f$-ideals of degree $2$.}

%The {\em total link vector} $l=l(\Delta)$ has components $l_i =
%\sum_{dim L_F = i}f(L_F)$. Where $f(\Delta)=\sum_{i=-1}^d f_i$ is
%known as the flag vector.

Consider a polynomial ring $S=k[x_1,\ldots,x_n]$ over a field $k$.
We say that a monomial ideal $I=(g_1,\ldots,g_m)\subset S$ where
$g_1,\ldots ,g_m$ are square-free monomials in $x_1,\ldots , x_n$
that form a minimal set of generators for $I$ is a {\em pure
square-free monomial ideal of degree $d$} if
$Supp(I)=\{x_1,\ldots,x_n\}$ and all the monomials $g_i\in S_d$ for
some $d>0$, where $S_d$ is the graded component of $S$, or in other
words all the $m_i$'s are of the same degree.
\begin{Definition}
{\em A square-free monomial ideal $I\subset S$ is said to be an {\em
$f$-ideal} if and only if both $\delta_{\mathcal{F}}(I)$ and
$\delta_{\mathcal{N}}(I)$ have the same $f$-vector. }
\end{Definition}
There is a natural question to ask:{\em  characterize all the
$f$-ideals in $S$}. Here, we precisely give the characterization of
$f$-ideals of degree $2$. So it is still open to characterize all
the $f$-ideals
of degree $\geq 3$.\\
%Before giving the proof of our main theorem, we need to define few
%notions. For a simplex $\Delta_{[n]}$ on the vertex set $[n]$:
%$$F_{d}(\Delta_{[n]}) =\big(\begin{array}{c}
 %          n \\
  %         d+1
   %      \end{array}\big)
 %$$
% where $F_d$ denotes the total $d$-dimensional faces in $\Delta_{[n]}$.
\begin{Lemma}\label{p1}
For a pure square-free monomial ideal $I=(g_1,\ldots,g_m)$ in
$S=k[x_1,\ldots,x_n]$ of degree $d$, the following equality holds;
$${n\choose d}=f_{d-1}(\delta_{\mathcal{F}}(I)) + f_{d-1}(\delta_{\mathcal{N}}(I))$$
\end{Lemma}
\begin{proof}
Let us take $I=(g_1,\ldots,g_m)\subset S$ be a square free monomial
ideal of degree $d$,  where $\{g_1,\ldots,g_m\}$ is the minimal set
of generators for $I$ and $deg(g_i)=d$ for all $i\in
\{1,\ldots,m\}$. So, corresponding to $I$ its facet simplicial
complex $\delta_{\mathcal{F}}(I)$ has
$$f_{d-1}(\delta_{\mathcal{F}}(I))= m.$$
As non-face complex $\delta_{\mathcal{N}}(I)$ will have the $d-1$
dimensional face $\{v_{i1} , \ldots , v_{id}\}$  if and only if
$x_{i1} \dots x_{id} \not\in I$ clear from the definition \ref{fc}.
So $\delta_{\mathcal{N}}(I)$ will have those $d-1$ dimensional faces
which are not appearing in $\delta_{\mathcal{F}}(I)$ because $I$ is
a pure square-free monomial ideal of degree $d$. Also, for a
simplicial complex on $n$ vertices the possible $d-1$ dimensional
faces are ${n\choose d}$. Therefore,
$f_{d-1}(\delta_{\mathcal{N}}(I)) = {n\choose d }-
f_{d-1}(\delta_{\mathcal{F}}(I))$.
\end{proof}
\begin{Remark}{\em
For instance, in example \ref{ex} one can see that $I=(xy,yz)\subset
k[x,y,z]$ is a pure square-free monomial ideal and;
$${3\choose 2}= f_{1}(\delta_{\mathcal{N}}(I)) +
f_{1}(\delta_{\mathcal{F}}(I))$$
$${3\choose 2} = 1 + 2.$$

%Note that the above equality holds for $f_{d-1}$ but not necessarily
%for all $f_i$ such that $i<d-1$. Moreover, we have taken $I$ to be
%pure square-free monomial ideal of degree $d$ which give rise a
%pure simplicial complex $\delta_{\mathcal{F}}(I)$.
 %One may think, what
%happened if $I$ will not necessarily be a pure monomial
%ideal?\\
%The next result gives an answer to the question: {\em For a given
%square-free monomial ideal $I\subset S$, when $dim
%(\delta_{\mathcal{F}}(I)) = dim (\delta_{\mathcal{N}}(I))$ }?
}
\end{Remark}
\begin{Lemma}\label{l1}
 Let $I$ be a square-free monomial ideal in $S$,
$$dim (\delta_{\mathcal{F}}(I)) = dim (\delta_{\mathcal{N}}(I))$$
 if and only if $ht(I)+deg(I)=n$.
\end{Lemma}
\begin{proof}
%The degree of the square-free monomial ideal $I\subset S $ gives us
%the dimension of the facet simplicial complex:
%$\delta_{\mathcal{F}}(I)$,\\
From definition \ref{fc}(a) it is clear that
$$dim(\delta_{\mathcal{F}}(I))= deg(I)-1.$$  Also, from Proposition 5.3.10 of \cite{Vi} it is clear that,
$$dim(\delta_{\mathcal{N}}(I))=n-ht(I)-1.$$
Which concludes the proof.
\end{proof}
One dimensional simplicial complexes on the vertex set $[n]$ are the
simple graphs. Also for a one dimensional simplicial complex the
ideal $I_{\mathcal{F}}$ is same as the edge ideal of a graph, for
details see \cite{Vi}.

%\begin{Corollary}\label{cor}
%Let $I\subset S$ be a square-free pure monomial ideal of degree
%$2$, $\delta_{\mathcal{F}}(I)$ is a connected graph if and only if
%$Supp(I)=\{x_1,\ldots,x_n\}$.
%\end{Corollary}
%\begin{proof} Let $I=(g_1,\ldots,g_m)\subset S$ be a square-free pure monomial ideal of degree
%$2$ with $Supp(I)=\{x_1,\ldots,x_n\}$.
%As there is one-to-one
%correspondence of  square-free monomial ideals and facet simplicial
%complexes.
% So $\delta_{\mathcal{F}}(I)$ will be one dimensional
%simplicial complex (which is a graph as well). Also for every vertex
%$v_i$ corresponding to $x_i$, there is a generator $g_t$ for some
%$t\in[n]$ such that $x_i/g_t$ (because
%$Supp(I)=\{x_1,\ldots,x_n\}$). Which implies
%$\delta_{\mathcal{F}}(I)$ is a connected graph. Conversely, suppose
%$\delta_{\mathcal{F}}(I)$ is a connected graph on $n$ vertices. As
%for any $x_j\in \{x_1,\ldots,x_n\}$ there is a vertex $v_j$ of
%$\delta_{\mathcal{F}}(I)$. Because $\delta_{\mathcal{F}}(I)$ is a
%connected graph, so there is an edge $\{v_j,v_k\}\in
%\delta_{\mathcal{F}}(I)$. This implies $x_jx_k=g_l$ for some $l\in
%[n]$. Hence $Supp(I)=\{x_1,\ldots,x_n\}$.
%\end{proof}
Our main theorem is as follows:
\begin{Theorem}\label{Th}
A pure square-free monomial ideal $I=(g_1,\ldots,g_m)\subset S$
of degree $2$  will be an {\em $f$-ideal}  if and only if:\\

($i$) $I$ is unmixed with $ht(I)= n-2$ ,\\

($ii$) ${n\choose 2}\equiv 0\, \, \, \, \, (\hbox{\em mod}\, 2)$  and\\

 ($iii$) $m=\frac{1}{2}{n\choose 2}$
\end{Theorem}
\begin{proof}
Suppose $I=(g_1,\ldots,g_m)\subset S$ is a pure square-free monomial
ideal  of degree $2$ and let $I$ be an $f$-ideal. So we have
$dim(\delta_{\mathcal{N}}(I))=1=dim(\delta_{\mathcal{F}}(I))$ which
by Lemma \ref{l1} implies $ht(I) = n- 2$. As $I$ is a pure
square-free monomial ideal of degree $2$, $\delta_{\mathcal{F}}(I)$
is a graph on the vertex set $[n]$ with no isolated vertex. So,
since $f(\delta_{\mathcal{F}}(I))=f(\delta_{\mathcal{N}}(I))$,
$\delta_{\mathcal{N}}(I)$ needs to be a graph on the same vertex set
$[n]$ with no isolated vertex. As $I_{\mathcal{N}}=\cap P_F$ where
the intersection is over all facets $F$ of $\delta_{\mathcal{N}}(I)$
by \cite{Vi} 5.3.10, which implies $I$ is unmixed of height $n-2$.

As $f_i$ denotes the number of $i$-dimensional faces, so
$$f_1(\delta_{\mathcal{F}}(I)) = m$$

 where $m$ is the number of monomial
generators of $I$. Also from Lemma \ref{p1}:
$$f_1(\delta_{\mathcal{N}}(I))= {n\choose 2} - m$$
As $I$ is an $f$-ideal, so we have ${n\choose 2} = 2m\equiv 0 \, \, \, \, \, (\hbox{\em mod}\, 2)$.\\
Conversely, let us take the pure square-free monomial ideal
$I=(g_1,\ldots,g_m)\subset S$ of degree $2$ satisfying the
conditions ($i$), ($ii$) and ($iii$). The simplicial complexes
$\delta_{\mathcal{F}}(I)$ and $\delta_{\mathcal{N}}(I)$ will have
the same $f$-vector follows immediately from Lemma \ref{p1} and
Corollary \ref{l1}.
\end{proof}
We conclude this paper with the following example.\\
\begin{Example}{\em Let $I=(x_1x_2,x_2x_3,x_3x_4)\subset
S=k[x_1,x_2,x_3,x_4]$ be a pure square-free monomial ideal
satisfying the conditions given in the above theorem, and
$\delta_{\mathcal{N}}(I)$ and $\delta_{\mathcal{F}}(I)$ are as
follows:

 \begin{figure}[h]
\begin{center}
  \includegraphics[width=8cm]{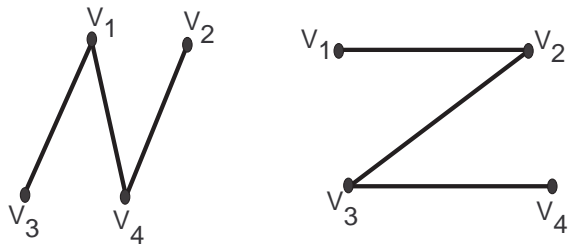}\vspace{0.0cm}
  \caption{Non-face and facet complex}\label{}
\end{center}
\end{figure}
Clearly, both the simplicial complexes have the same $f$-vector
$(4,3)$. Hence $I$ is an $f$-ideal of degree $2$.

}
\end{Example}


\vspace{1 pt}
\begin{thebibliography}{1}
%\bibitem{AP} S. Ahmad,\ D. Popescu, Sequentially Cohen-Macaulay monomial ideals of embedding dimension
%four, Bull. Math. Soc. Sc. Math. Roumanie, 50(98),no 2,
%(2007),99-110( see www.rms.unibuc.ro/bulletin or
%Arxiv:Math.AC/0702569).

%\bibitem{AP1}I. Anwar, D. Popescu, Stanley Conjecture in small embedding
%dimension, J.Algebra {\bf{318}} (2007) $1027 - 1031$

%\bibitem{Ap} J. Apel,  On a conjecture of R.P. Stanley, in
%{\em Part I - Monomial
%ideals.J.Algebra.Comb.},{\bf{17}},36-59(2003).

\bibitem{BH} W. Bruns, J. Herzog, {\em Cohen Macaulay rings},
Vol.39, Cambridge studies in advanced mathematics, revised
edition,1998.


\bibitem{F1} S. Faridi, {\em The facet ideal of a simplicial complex}, Manuscripta Mathematica,
109, (2002), 159-174.

\bibitem{F2} S. Faridi, {\em Simplicial Tree are sequentially Cohen-Macaulay}, Arxiv:Math.AC/0702569.

%\bibitem{Fi} J. Fine, {\em A complete $h$-vector for convex polytopes }, Arxiv:Math.AC/5722v1.


%\bibitem{Gr}W. Gr\"{o}bner, \"{U}ber, die algebraischen
%Eigenschaften  der Integrale von linearen Differentialgleichungen
%mit konstanten Koeffizienten. Monatsh. Math. Phys. {\bf{47}},
%247-284(1939)(German)

%\bibitem{Gr1} W. Gr\"{o}bner, \"{U}ber, die Eliminationstheorie.
%Monatsh. Math.{\bf{54}}, 71-78(1950)(German)


%\bibitem{HJY} J. Herzog,\ A. Soleyman Jahan,\ S. Yassemi, Stanley
%decompositions and partitionable simplicial complexes, to appear in J. Algebraic Cobinatorics, Arxiv:Math. AC/0612848v2.

%\bibitem{HP} J. Herzog,\ D. Popescu, Finite filtrations of modules
%and shellable multicomplexes, Manuscripta Math. {\bf 121}, (2006),
%385-410.

%\bibitem{HT} S.Hosten, R.R.Thomas, Standard pairs and group
%relaxations in integer programming, J.Pure Appl. Alg.
%{\bf{139}}, (1999), 133-157.

\bibitem{Mi}E. Miller, B. Sturmfels, 'Combinatorial Commutative Algebra',
Springer-Verlag New York Inc. 2005.

%\bibitem{MS} D. Maclagan,\ G. Smith,  Uniform bounds on multigraded
%regularity J.Alg.Geom. \textbf{14}(2005), 137-164.


%\bibitem{Pl} W.Plesken, {\it{Janet's Algorithm.}} In: S. Abenda, G.
%Gaeta, S. Walcher (eds), Symmetry and perturbation Theory (SPT
%2002). Singapore 2003.

%\bibitem{Po} D. Popesscu, {\it{Stanley depth of mononmial
%ideals}}, to appear in the precceedings of 6th Congress of Romanian
%Mathematicians, 2007.

%\bibitem{PR} W.Plesken, D.Robertz, {\it{Janet's approach to presentations
%and resolutions for polynomials and linear pdes.}}, Arch. Math.
%84(2005) 22-37.

%\bibitem{Ra} A. Rauf, Stanley decompositions, pretty clean filtrations and reductions modulo regular
%elements, Bull. Math. Soc. Sc. Math. Roumanie, 50(98),no 4, (2007) (see
%www.rms.unibuc.ro/bulletin).

\bibitem{St} R. P. Stanley, {\em Combinatorics and commutative
algebra}, Second edition. Progress in Mathematics, 41. Birkhuser
Boston, MA, 1996, x+164 pp. ISBN: 0-8176-3836-9.

\bibitem{St1} R. P. Stanley,{\em Cohen-Macaulay Rings and constructible
polytopes}, Bull. Amer. Math. Soc. 81(1975),133-142.

\bibitem{Vi} R. H. Villarreal, {\em Monomial algebras}, Dekker, New York,
2001.


\end{thebibliography}
\end{document}